\input amstex\documentstyle{amsppt}  
\pagewidth{12.5cm}\pageheight{19cm}\magnification\magstep1
\topmatter
\title Exceptional representations of Weyl groups\endtitle
\author G. Lusztig\endauthor
\address{Department of Mathematics, M.I.T., Cambridge, MA 02139}\endaddress
\thanks{Supported in part by National Science Foundation grant DMS-1303060.}\endthanks
\endtopmatter   
\document

\define\Irr{\text{\rm Irr}}

\define\da{\dagger}

\define\si{\sim}

\define\bX{\bar X}

\define\op{\oplus}
   
\define\part{\partial}
\define\emp{\emptyset}

\define\n{\notin}
\define\iy{\infty}
\define\m{\mapsto}
\define\do{\dots}

\define\sub{\subset}    

\define\T{\times}
\define\ti{\tilde}
\define\nl{\newline}
\redefine\i{^{-1}}

\define\ot{\otimes}

\define\sg{\text{\rm sgn}}
\define\tr{\text{\rm tr}}

\define\g{\gamma}
\redefine\d{\delta}
\define\e{\epsilon}

\define\r{\rho}

\define\x{\xi}

\redefine\G{\Gamma}

\define\Ph{\Phi}

\define\boc{\bold c}

\define\kk{\bold k}

\define\CC{\bold C}

\define\FF{\bold F}

\define\NN{\bold N}

\define\RR{\bold R}

\define\ZZ{\bold Z}

\define\ca{\Cal A}
\define\cb{\Cal B}

\define\cd{\Cal D}
\define\ce{\Cal E}

\define\co{\Cal O}

\define\cu{\Cal U}

\define\cx{\Cal X}

\define\tn{\ti n}

\define\bS{\bar S}

\define\AL{AL}
\define\BC{BC}
\define\BL{BL}
\define\DL{DL}
\define\KL{KL}
\define\SPE{L1}
\define\OBC{L2}
\define\ORA{L3}
\define\LCE{L4}
\define\HEC{L5}
\define\BAR{L6}
\define\LV{LV}

\subhead 1.1\endsubhead
Let $W$ be a finite, irreducible Coxeter group and let S be the set of simple
reflections of $W$; let $l:W@>>>\NN$ be the length function.
Let $\Irr W$ be a set of representatives for the 
isomorphism classes of irreducible representations of $W$ over $\CC$, the 
complex numbers. Let $\ca=\CC[v,v\i]$, $\ca'=\CC[v^2,v^{-2}]$ ($v$ an indeterminate). 
We have $\ca'\sub\ca\sub K\supset K'\supset\ca'$ where $K=\CC(v),K'=\CC(v^2)$. 
Let $H$ be the Hecke algebra over $\ca$ associated to $W$; thus 
$H$ has generators $T_s (s\in S)$ and generators $(T_s+1)(T_s-v^2)=0$ for 
$s\in S$, $T_sT_{s'}T_s\do=T_{s'}T_sT_{s'}\do$ for $s\ne s'$ in $S$ (both 
products have $m$ factors where $m$ is the order of $ss'$ in $W$). 
Let $H'$ be the $\ca'$-subalgebra of $H$ generated by 
$T_s (s\in S)$; note that $H=\ca\ot_{\ca'}H'$. 
Let $H_K=K\ot_\ca H$, $H_{K'}=K'\ot_{\ca'}H'$ so that $H_K=K\ot_{K'}H_{K'}$.
It is known \cite{\OBC} (see also 1.2 below) that the algebra $H_K$ is canonically 
isomorphic 
to the group algebra $K[W]$. Hence any $E\in\Irr W$ can be viewed as a simple
$H_K$-module $E_v$. We say that $E$ is {\it ordinary} if $E_v$ is obtained by 
extension of scalars from an $H_{K'}$-module; otherwise, we say that $E$ is
{\it exceptional}. Let $\Irr_0W$ (resp. $\Irr_1W$) be the set of all 
$E\in\Irr W$ which are ordinary (resp. exceptional).

We define a subset $\ce W$ of $\Irr W$ as follows.
If $W$ is not of type $E_7,E_8,H_3,H_4$, we set $\ce W=\emp$.
If $W$ is of type $E_7,E_8,H_3,H_4$, then $\ce W$ consists of 
$2^a$ representations of dimension $2^b$ where $2^a=2$ for $E_7,H_3$, $2^a=4$ 
for $E_8,H_4$ and $2^{a+b}$ is the largest power of $2$ that divides the order
of $W$; thus $2^b$ is $512,4096,4,16$ respectively. 

When $W$ is crystallographic we have $\Irr W-\ce W\sub\Irr_0W$ (see \cite{\BC}) and $\ce W\sub\Irr_1W$ (a
result of Springer); hence $\Irr W-\ce W=\Irr_0W$ and $\ce W=\Irr_1W$. The same holds when $W$ is not 
crystallographic. (The fact $\ce\sub\Irr_1W$ for $W$ of type $H_3$ was pointed out in \cite{\OBC}.
The fact that any $E\in\Irr W-\ce$ is ordinary for $W$ of type $H_4$ can be seen from the fact that,
according to \cite{\AL}, $E$ can be realized by a $W$-graph which is even (in the sense 
that the vertices can be partitioned into two subsets so that no edge connects vertices in the same subset).

In this paper we try to understand various consequences in representation 
theory of the existence of exceptional representations.

\subhead 1.2\endsubhead
Let $\{c_w;w\in W\}$ be the basis of $H$ which in \cite{\KL} was denoted by
$\{C'_w;w\in W\}$. Let $\le_{LR},\le_L$ be the preorders on $W$ defined in 
\cite{\KL} and let $\si_{LR},\si_L$ be the corresponding equivalence relations
 on $W$
(the equivalence classes are called the two-sided cells and left cells
respectively). For $x,y\in W$ we write $c_xc_y=\sum_{z\in W}h_{x,y,z}c_z$.
For $z\in W$ there is a unique number $a(z)\in\NN$ such that for any $x,y$ in 
$W$ we have $h_{x,y,z}=\g_{x,y,z\i}v^{a(z)}\mod v^{a(z)-1}\ZZ[v\i]$
where $\g_{x,y,z\i}\in\NN$ and $\g_{x,y,z\i}>0$ for some $x,y$ in $W$.
Moreover, $z\m a(z)$ is constant on any two-sided cell. (See \cite{\HEC}.)
Let $J$ be the $\CC$-vector space with basis $\{t_w;w\in W\}$. It has an
associative $\CC$-algebra structure given by 
$t_xt_y=\sum_{z\in W}\g_{x,y,z\i}t_z$; it has a unit element of the form
$\sum_{d\in\cd}t_d$ where $\cd$ is a subset of $W$ consisting of certain 
involutions (that is elements with square $1$). (See \cite{\HEC}.) 
Let $h\m h^\da$ be the algebra automorphism of $H$ such that $T_s^\da=-T_s\i$ for $s\in S$.
Now the $\ca$-linear map $H@>>>\ca\ot J$ given by
$c_x^\da\m\sum_{d\in\cd,z\in W,d\si_Lz}h_{x,d,z}t_z$ induces an algebra 
isomorphism $H_K@>>>K\ot J$ and (by specializing $v=1$) an algebra
isomorphism $\CC[W]@>>>J$ hence an algebra isomorphism $K[W]@>>>K\ot J$. (See \cite{\HEC}.)
Now if $E\in\Irr W$ then $E_v$ in 1.1 is obtained as follows.
We first view $K\ot E$ as a $K\ot J$-module $E_\iy$ via the isomorphism
$K[W]@>>>K\ot J$ above and then view $E_\iy$ as an $H_K$-module $E_v$ via the
isomorphism $H_K@>>>K\ot J$. Note that for $x\in W$ we have
$$\tr(c_x^\da,E_v)=\sum_{d\in\cd,z\in W;d\si_Lz}h_{x,d,z}\tr(t_z,E_\iy).\tag a$$
We show:

(b) {\it If $x\in W$ satisfies $x^2=1$, or more generally, if $x\si_L x\i$ then
there exists $E\in\Irr W$ such that $\tr(t_x,E_\iy)\ne0$.}
\nl
It is enough to show that $\sum_{E\in\Irr W}\tr(t_x,E_\iy)\tr(t_{x\i},E_\iy)\ne0$. The last sum
is equal to the trace of the $K$-linear map $K\ot J@>>>K\ot J$, 
$\x\m t_x\x t_{x\i}$ (we use that $t_w\m t_{w\i}$ defines an isomorphism of the algebra $K\ot J$
onto the algebra with opposed multiplication) hence it is equal to
$\sum_{u,u'\in W}\g_{x,u,u'{}\i}\g_{u',x\i,u\i}$. Thus it is enough to show that
the last sum is $\ne0$. Now each term in the last sum is in $\NN$ hence it is
enough to show that for some $u,u'$ we have
$\g_{x,u,u'{}\i}\g_{u',x\i,u\i}>0$. We take $u=x\i$ and $u'=d$ where
$d$ is the unique involution in $\cx$ such that $x\si_L d$. It is enough to show that 
$\g_{x,x\i,d}\g_{d,x\i,x}>0$. But the last product is $1$ since $x\si_Ld\si_Lx\i$ so that
$\g_{x,x\i,d}=1,\g_{d,x\i,x}=\g_{x\i,x,d}=1$. (See \cite{\HEC}.) This proves (b).

\subhead 1.3\endsubhead
If $E\in\Irr W$ then there is a unique two-sided cell $c$ such that
$t_x:E_\iy@>>>E_{\iy}$ is nonzero for some $x\in c$. This gives us a (surjective) map
$E\m c$ from $\Irr W$ to the set of two-sided cells; its fibre at a two-sided cell $c$
is denoted by $\Irr^cW$. One checks that if some $E\in\Irr^cW$ is exceptional then any
$E\in\Irr^cW$ is exceptional; in this case we say that $c$ is exceptional. If some/any
$E\in\Irr^cW$ is ordinary, we say that $c$ is ordinary. 
An involution $x$ in $W$ is said to be ordinary (resp. exceptional) if $l(x)=a(x)\mod2$
(resp. $l(x)=a(x)+1\mod2$). Note that any two-sided cell $c$ contains some ordinary involution
(for example, $c\cap\cd$ is a nonempty set consisting of ordinary involutions). We show:

(a) {\it If $c$ is an ordinary two-sided cell, then for any $x\in c$ such that $x\si_Lx\i$
we have $l(x)=a(x)\mod2$. In particular, any involution in $c$ is ordinary.}
\nl
By 1.2(a) we can find $E\in\Irr W$ such that $\tr(t_x,E_\iy)\ne0$. By definition we have $E\in\Irr^cW$
hence $E$ is ordinary; since $v^{l(x)}c_x^\da\in H'$, it follows that $\tr(v^{l(x)}c_x^\da,E_v)\in\CC(v^2)$.
Using this and 1.2(a) we deduce
$$v^{l(x)}\sum_{d\in\cd,z\in W;d\si_Lz}h_{x,d,z}\tr(t_z,E_\iy)\in\CC(v^2).\tag b$$
Let $a_0$ be the value of the $a$-function on $c$. For $z,d$ in the last sum such that
$\tr(t_z,E_\iy)\ne0$ we have $z\in c$ hence $a(z)=a_0$ and $h_{x,d,z}=\g_{x,d,z\i}v^{a_0}$ plus a 
$\ZZ$-linear combination of strictly smaller powers of $v$; moreover we have $\g_{x,d,z\i}=\g_{z\i,x,d}$
and this is $1$ if $z=x$ and $d$ is the unique element of $\cd$ such that $d\si_Lx$ and is $0$ otherwise.
Thus (b) becomes
$$v^{l(x)}v^{a_0}\tr(t_x,E_\iy)+\text{lin.comb.of strictly smaller powers of }v\in\CC(v^2).$$
Since $\tr(t_x,E_\iy)\in\CC-\{0\}$ it follows that $l(x)+a_0\in2\ZZ$ and (a) follows.

We now show:

(c) {\it If $c$ is an exceptional two-sided cell, then $c$ contains both ordinary and exceptional
involutions. More precisely, any left cell in $c$ contains exactly one ordinary involution and
exactly one exceptional involution.}
\nl
Let $n_c$ (resp. $\tn_c$) be the number of ordinary (resp. exceptional) involutions in $c$. 
Let $c'=w_0c$ where $w_0$ is the longest element of $W$. Then $c'$ is again an exceptional two-sided cell.
In type $E_7$ or $H_3$ we have $c'=c$. Since $w_0$ is central in $W$ and of odd length, for any involution 
$x$ in $c$, $w_0x$ is again an involution in $c$ and $x$ is ordinary if and only if $w_0x$ is exceptional; 
thus we have $n_c=\tn_c$. In type $E_8$ or $H_4$ we have $c'\ne c$; more precisely the value of the 
$a$-function on $c$ has a different parity than that on $c'$. Since $w_0$ is central in $W$ and of even 
length, for any involution $x$ in $c$, $w_0x$ is an involution in $c'$ and $x$ is ordinary if and only if 
$w_0x$ is exceptional; thus we have $n_c=\tn_{c'}$ and $\tn_c=n_{c'}$.

Note that $\Irr^cW$ consists of two elements of dimension $m_c$ where $m_c$ is the number of left cells in 
$c$. It is known that if $\G$ is any left cell in $c$ then $\G$ carries a $W$-module structure isomorphic to
the direct sum of the two representations in $\Irr^cW$. In type $E_7,E_8$, using \cite{\ORA, 12.15}, we 
deduce that $\G\cap\G\i$ has exactly two elements (a similar result can be proved in type $H_3,H_4$). The 
unique element of $\G\cap\cd$ is one of these two elements and is an ordinary involution. Also
any involution in $\G$ is contained in $\G\cap\G\i$. We see that $n_c\ge m_c\ge\tn_c$. 
In type $E_7,H_3$ we have $n_c=\tn_c$ hence $n_c=m_c=\tn_c$; we see that (c) holds in this case.
In type $E_8,H_4$ we have $n_c\ge\tn_c$ (and similarly $n_{c'}\ge\tn_{c'}$).
Using $n_c=\tn_{c'}$ and $\tn_c=n_{c'}$ we deduce $n_c=\tn_{c'}=\tn_c=n_{c'}=m_c=n_{c'}$; we see that
(c) holds in this case. 

\subhead 1.4\endsubhead
In this subsection we assume that $W$ is crystallographic.
Let $G$ be a simple algebraic group over an algebraic closure $\kk$ of a finite field $\FF_q$ with $q$
elements with a fixed split $\FF_q$-structure such that the Weyl group of $G$ is $W$ in 1.1. The variety 
$\cb$ of Borel subgroups of $G$ has a natural $\FF_q$-structure with Frobenius map $F:\cb@>>>\cb$. For each 
$w\in W$ let $\co_w$ be the $G$-orbit on $\cb\T\cb$ (diagonal action) indexed by $w$ and let $\bX_w$ be 
the closure in $\cb$ of the variety $\{B\in\cb;(B,F(B))\in\co_w\}$ of \cite{\DL}. Now $G(\FF_q)$ acts 
naturally on the
$l$-adic intersection cohomology spaces $IH^i(\bX_w)$. An irreducible representation of $G(\FF_q)$ is said to
be unipotent if it appears in the $G(\FF_q)$-module $IH^i(\bX_w)$ for some $w,i$.
Let $\cu_q$ be the a set of representatives for the isomorphism classes of unipotent representations of
$G(\FF_q)$. Let $\r\in\cu_q$. By \cite{\ORA, 3.8}, for any $\r\in\cu_q$, any $z\in W$ and any $j\in\ZZ$ we 
have
$$\align&(\r:IH^j(\bX_z))_{G(\FF_q)}\\&=\text{coefficient of }v^j\text{ in }
(-1)^j\sum_{E\in\Irr W}c_{\r,E}\tr(v^{l(z)}c_z,E(v))\endalign$$
where $c_{\r,E}$ are uniquely defined rational numbers and $\tr(v^{l(z)}c_z,E(v))\in\ca$.
Moreover, by \cite{\ORA, 6.17}, given $\r$ as above, there is a unique 
two-sided cell $\boc_\r$ of $W$ such that $c_{\r,E}=0$ whenever $E\n\Irr^{\boc_\r}W$.
For a two-sided cell $\boc$ we write $\cu^\boc_q=\{\r\in\cu_q;\boc_\r=\boc\}$. We see that for 
any $\r\in\cu_q^\boc$, any $z\in W$ and any $j\in\ZZ$ we have
$$\align&(\r:IH^j(\bX_z))_{G(\FF_q)}\\&=\text{coefficient of }v^j\text{ in }
(-1)^j\sum_{E\in\Irr^\boc W}c_{\r,E}\tr(v^{l(z)}c_z,E(v)).\tag a\endalign$$
To any $\r\in\cu_q$ we associate a sign $\e_\r\in\{1,-1\}$ by the following requirement:
if $\r$ appears in $IH^j(\bX_z)$ with $z\in W,j\in\ZZ$ then $\e_\r=(-1)^j$; this is well defined by 
\cite{\ORA, 6.6}.
We say that $\r\in\cu_q$ is ordinary (resp. exceptional) if $\e_\r=1$ (resp. $\e_\r=-1$).
We show:

(b) {\it If $\boc$ is an ordinary two-sided cell then any $\r\in\cu_q^\boc$ is ordinary. If $\boc$ is an 
exceptional two-sided cell, then $\cu_q^\boc$ consists of two ordinary and two exceptional representations.}
\nl
Assume first that $\boc$ is ordinary. Since $\tr(v^{l(z)}c_z,E(v))\in\ca'$ for $E\in\Irr^\boc W$, $z\in W$, 
we see from (a) that for any $\r\in\cu^\boc_q$ we have $(\r:IH^j(\bX_z))_{G(\FF_q)}=0$ if $j$ is odd. Thus 
$\r$ is ordinary. Assume next that $\boc$ is exceptional. Then $\cu^\boc_q$ consists of four representations
of which two appear in $IH^0(\bX_1)$ hence are ordinary and the other two appear in $IH^7(\bX_z)$
where $z$ is an element of length $7$ in $W$.

\subhead 1.5\endsubhead
Let $S_i$ be the $i$-th symmetric power of the reflection representation of $W$ and let $S=\op_iS_i$, a
commutative algebra over $\RR$. Let $I$ be the ideal of $S$ generated by the $W$-invariant elements of $S$
of degree $>0$. Let $\bS=S/I$ and let $\bS_i$ be the image of $S_i$ in $\bS$. Note that $\bS_i$ is a 
$W$-module. For any $E\in\Irr W$ we set $P_E=\sum_{i\ge0}(E:\bS_i)X^i\in\NN[X]$. We note the following 
property:

(a) {\it If $E$ is ordinary then $P_E$ is palindromic. If $E$ is exceptional then $P_E$ is not palindromic.} 
\nl
(A polynomial $P(X)\in\CC[X]$ is said to be palindromic if there exists
$u\in\NN$ such that $P(X\i)=X^{-u}P(X)$.) When $W$ is crystallographic this has been noted in
\cite{\BL}. When $W$ is dihedral or of type $H_3$ this is easily verified. 
When $W$ is of type $H_4$ this follows from \cite{\AL}. We will now give an explanation for why (a) holds
assuming that $W$ is crystallographic. 

Let $\boc$ be the two-sided cell such that $E\in\Irr^\boc W$. 
It is known that $\cu_q^\boc$ (see 1.4) can be naturally indexed by a set independent
of $q$ so that when $\r\in\cu_q^\boc$, the dimension of $q$ can be regarded as a polynomial $\d_\r$ in $q$ 
with rational coefficients; more precisely, we have $d_\r(X)=e\i X^{a_0}(X-1)^sf(X)$ where $e\in\ZZ_{>0}$,
$a_0\in\NN$ depends only on $\boc$, not on $\r$, $s\in\NN$ is such that $\e_\r=(-1)^s$ and $f$ is a 
product
of cyclotomic polynomials $\Ph_r(X)$ with $r\ge2$. Also the degree of the polynomial $X^{a_0}(X-1)^sf(q)$ is
a number $A_0$ depending only on $\boc$, not on $\r$. It follows that $d_\r(X\i)=(-1)^sX^{-a_0-A_0}d_\r(X)$.

We now assume that $\boc$ is ordinary. Then we have $\e_\r=1$ for each $\r$ as above (hence $s$ is even), 
see 1.4(b). From \cite{\ORA, 4.23} it is known that $P_E(X)$ is a constant linear combination of polynomials
$d_\r(X)$ with $\r\in\cu^\boc_q$. Since $d_\r(X\i)=X^{-a_0-A_0}d_\r(X)$ for each $\r$ it follows that
$P_E(X\i)=X^{-a_0-A_0}P_E(X)$.

\subhead 1.6\endsubhead
In this subsection we assume that $W$ is crystallographic.
Let $\boc$ be a two-sided cell. Let $E_\boc$ be the special representation in $\Irr^\boc W$ (see 
\cite{\ORA}). For each left cell $\G$ in $\boc$ let $[\G]$ be the
$W$-module carried by $\G$. Let $[[\boc]]$ be the $W$-module carried by the set of involutions in $\boc$
defined in \cite{\LV}. We have the following result.

(a) {\it Assume that $\boc$ is ordinary. There is a unique $E\in\Irr^\boc W$ such that $E$ appears in $[\G]$
for every $\G$ as above and $E$ appears in $[[\boc]]$, namely $E=E_\boc$.}
\nl
If $W$ is of classical type, then it is known that $[[\boc]]$ is a sum of copies of $E_\boc$ and that
$E_\boc$ appears in each $[\G]$ with multiplicity one. Hence (a) holds in this case. We now assume that $W$ 
is of exceptional type. If $\boc$ is not the two-sided cell containing $E$ of dimension
$4480$ (in $E_8$) or $12$ (in $F_4$) or $2$ (in $G_2$) then there is exactly one $E$ which appears in
each $[\G]$ namely $E_\boc$ and $E_\boc$ appears in $[[\boc]]$, see \cite{\LCE}; hence (a) holds in this 
case. 
We now assume that $\boc$ is the two-sided cell containing $E$ of dimension
$4480$ (in $E_8$) or $12$ (in $F_4$) or $2$ (in $G_2$).
Then there are exactly two $E$ which appear in each $[\G]$ namely $E_\boc$ and the $E$ of dimension
$7168$ (in $E_8$) or $16$ (in $F_4$) or $2$ (non-special) in $G_2$, see \cite{\LCE}. Now one verifies that 
$E_\boc$ appears
in $[[\boc]]$ but $7168$ (in $E_8$) or $16$ (in $F_4$) or $2$ (non-special) in $G_2$ do not appear in 
$[[\boc]]$. Hence (a) holds in this case. 

Note that if $\boc$ is exceptional then there are exactly two $E\in\Irr^\boc W$ such that $E$ appears in 
$[\G]$ for every $\G$ as above and $E$ appears in $[[\boc]]$; one of them is $E=E_\boc$.

Let $\sg$ be the sign representation of $W$. The following result has been noted in \cite{\SPE}.

(b) {\it If $\boc$ is ordinary then $E_\boc\ot\sg$ is a special representation.
If $\boc$ is exceptional then $E_\boc\ot\sg$ is not a special representation.}
\nl
Note that the first statement of (b) can be deduced from (a) applied to $w_0\boc$ (an ordinary two-sided
cell) since if $\G'$ is a left cell in $w_0\boc$ then $[\G']\cong[\G]\ot\sg$ for some left cell in $\boc$
and $[[w_0\boc]]\cong[[\boc]\ot\sg$ (this follows from the inversion formula in \cite{\BAR}).

\enddocument